\documentclass[12pt,a4paper]{article} %[g-shadowing-cor 03.12.22]   Updated
\usepackage{amsfonts,amssymb} 
\usepackage[cp1251]{inputenc}

\usepackage[margin=25mm, a4paper]{geometry}  % journal quality

\usepackage{tikz}  \def\bwait#1#2{#2}

\def\n{\noindent}  \def\?#1{}
\def\IZ{{\mathbb{Z}}}  \def\IR{{\mathbb{R}}}

   \def\Tor{{\mathbb{T}}}   \def\la{\lambda}  \def\cN{{\cal N}}

   \def\phi{\varphi}  \def\diam{{\rm diam}}
\def\ep{\varepsilon}  \def\t{\tilde} \def\bdelete#1{}
    \def\v#1{\vec{#1}}
\def\mod1{\,({\rm mod\ } 1)\,}

\def\beq#1#2{\begin{equation} \label{#1} #2 \end{equation}}
 
\def\toas#1{\stackrel{#1}{\longrightarrow}}
\def\function#1{\left\{\!\!\!\begin{array}{ll} #1 \end{array} \right.}
\def\proof{\smallskip \noindent {\bf Proof. \ }}       %start of proof
\def\blanksquare{\,\,\,$\sqcup\!\!\!\!\sqcap$}         %blank  square
\def\qed{\hfill\blanksquare\linebreak\smallskip\par}   %end of proof
\def\thname{Theorem}  \def\lmname{Lemma}    \def\prname{Proposition}
\def\dfname{Definition}  \def\crname{Corollary}  \def\rmname{Remark}
\def\exname{Example}  

\newtheorem{theorem}{\thname}[section]   %Numbering: Theorem--Other section
\newtheorem{lemma}{\lmname}[section]     %{lemma}[theorem]{Lemma}   subsection
\newtheorem{proposition}[lemma]{\prname} %lemma
\newtheorem{example}{\exname}[section]

\newtheorem{dftn}{\dfname}[section]
\newenvironment{definition}{\begin{dftn}\rm}{\end{dftn}} %section
\def\bdef#1{\begin{definition} #1 \end{definition}}
\newtheorem{rmrk}[lemma]{\rmname}
\newenvironment{remark}{\begin{rmrk}\rm}{\end{rmrk}}     %lemma

\makeatletter\def\fps@figure{htbp}\makeatother %figure pos: tbp - standard

\def\bcb#1{\textcolor{black}{#1}}           % blue
\begin{document}

\title{Average shadowing and gluing property}
\author{Michael Blank\thanks{
        Institute for Information Transmission Problems RAS
        (Kharkevich Institute);}
        \thanks{National Research University ``Higher School of Economics'';
        e-mail: blank@iitp.ru}
       }
\date{December 3, 2022} %{\today}%{February 27, 2022} %
\maketitle

\begin{abstract} The purpose of this work is threefold: 
(i) extend shadowing theory for discontinuous and non-invertible systems,
(ii) consider more general classes of perturbations (for example, small only on average), 
(iii) establish a general theory based on the property that the shadowing holds 
for the case of a single perturbation. The ``gluing'' construction used in the analysis 
of the last property turns out to be the key point of this theory. \end{abstract}

{\small\n
}
%http://aps.ecnu.edu.cn/UserFiles/File/MSC2020-Mathematica l% 37Bxx Topological dynamics
% 37B65 Approximate trajectories, pseudotrajectories, shadowing
% 37B05 Dynamical systems involving transformations and group actions with special properties (minimality, distality, proximality, expansivity, etc.)
% 37B10 Symbolic dynamics
% 37C50 Approximate trajectories (pseudotrajectories, shadowing, etc.) in smooth dynamics
%%%%%%%%%%%%%%%%%%%%%%%%%%%%%%%%

\section{ Introduction}

Due to the unavoidable presence of various errors and perturbations in the modeling of dynamical 
systems (DS), the question arises about the relationship between the asymptotic properties of the 
simulated system and the simulation results. This question is especially important (and complicated) 
in the case of chaotic dynamics. 
At the level of connections between individual trajectories of a hyperbolic system and the 
corresponding pseudo-trajectories\footnote{Approximate trajectories of a system under small 
     perturbations already considered by G. Birkhoff \cite[(1927)]{Bi} with a completely different purpose.}, 
this problem was first posed by D.V. Anosov \cite[(1967-70)]{An,An2} as a key step of the analysis 
of structural stability of diffeomorphisms. 
A similar but much less intuitive approach called ``specification'' in the same setting 
was proposed by R. Bowen \cite[(1975)]{Bo}. 
Informally, both approaches ensure that errors do not accumulate during the process of modeling: 
in the systems with shadowing property each approximate trajectory can be uniformly traced 
by a true trajectory on the arbitrary long period of time. Naturally, this is of great
importance in chaotic systems, where even an arbitrary small error in the starting position lead to 
(exponentially in time) large divergence of trajectories.

Further development (see main results, generalizations and numerous references in two 
monographes dedicated to this subject \cite{Pi, PS}) demonstrated deep relations between 
the shadowing property and various ergodic characteristics of dynamical systems. 
In particular, it has been shown that for a diffeomorphism the shadowing property implies 
the uniform hyperbolicity. 
To some extent, this limits the theory of uniform shadowing to an important but very special 
class of hyperbolic dynamical systems.

It is worth noting that a great body of results related to shadowing on finite time 
intervals (see \cite{Pi} and references therein), being of undeniable interest from an 
applied point of view, gives practically no information about (the most important) 
asymptotic as time goes to infinity properties of a chaotic system. 

To get out of this impasse, M.~Blank proposed  a new concept of average shadowing 
\cite[(1988)]{Bl}, which in particular takes into 
account much more general classes of perturbations (for example, of Gaussian type). 
Further generalizations of this concept and recent progress in this direction see   
\cite{LS, KO, KKO, WOC} (and an extensive list of other references therein). 
Potentially the average shadowing was expected to be applicable to non-hyperbolic 
systems, but this has not yet been proven. 
Moreover, K. Sakai \cite{Sa,Sa2} showed that a diffeomorphism satisfies the (one-sided)
average shadowing property only if it is positively expansive.\footnote{This property is  
     close to hyperbolicity (see discussion in Section~\ref{s:pre}). Here the one-sided 
     shadowing means that only forward semi-trajectories are taken into account.}

Our goal is to overcome these limitations in order to extend the theory of shadowing for 
discontinuous and non-invertible systems, as well as to consider more general classes 
of perturbations. To this regard, one of the main difficulties in the shadowing theory is the 
presence of infinitely many moments of perturbations which makes the analysis 
highly non-local. It is therefore very desirable to reduce the shadowing problem to the 
situation with a single perturbation, albeit with tighter control of the approximation accuracy. 
To carry out this program, we have developed a fundamentally new ``gluing'' construction, 
consisting in the approximation of segments of true trajectories.  See exact definitions 
and details of the construction in Section~\ref{s:pre}.

In what follows we will demonstrate that the ``gluing'' construction is applicable in 
all situations under study. Moreover, using it we are able to study various combinations 
of types of perturbations and types of shadowing.\footnote{Previously, a separate method 
    was developed to analyze each specific combination of a perturbation and a type of shadowing.}

The paper is organized as follows. In Section~\ref{s:pre} we give general definitions 
related to the shadowing property and introduce the key tool of our analysis -- the gluing 
property. In Section~\ref{s:main} we formulate and prove the main result -- Theorem~\ref{t:main}, 
which deduces various versions of shadowing from the gluing property. 
Finally, Section~\ref{s:ver} is devoted to the verification of the gluing property for various 
classes of discrete time dynamical systems, starting with hyperbolic and piecewise 
hyperbolic and non-invertible mappings and ending with systems with neutral singularities. 
The latter case demonstrates the difference between strong (\ref{e:glu}) and weak (\ref{e:glu-w})  
versions of the gluing property, which shows that even with uniformly small perturbations it is 
possible that only the average shadowing takes place (but not the uniform one).

\section{Preliminaries}\label{s:pre}
Let $T:X \to X$ be a non-necessarily invertible map from a metric space $(X,\rho)$ into itself. 

\bdef{A {\em trajectory} of the map $T$ starting at a point $x\in X$ is a sequence 
of points $\v{x}:=\{\dots,x_{-2},x_{-1},x_0,x_1,x_2,\dots\}\subset X$, 
for which $x_0=x$ and $Tx_i=x_{i+1}$ 
for all available indices $i$. The part of $\v{x}$ corresponding to non-negative indices 
is called the {\em forward (semi-)trajectory}, while the part corresponding to non-positive indices 
is called the {\em backward (semi-)trajectory}.}

An important comment here is that although the forward trajectory is always uniquely determined 
by $x=x_0$ and infinite, while the backward trajectory might be finite (if its ``last'' point has no 
preimages)\footnote{In this case we  are speaking only about available indices $i$ in the definition.} 
and that for a given $x=x_0$ there might be arbitrary many admissible backward trajectories. 

\begin{remark}
The reason for introducing the somewhat unusual notion of the backward trajectory 
of a non-invertible dynamical system is that when analyzing the connections between the 
trajectories of the original and perturbed systems, we inevitably have to go back and through 
in time. Therefore it is more convenient  to define infinite in both directions trajectories from 
the very beginning.
\end{remark}

\bdef{A {\em pseudo-trajectory} of the map $T$ is a sequence 
of points $\v{y}:=\{\dots,y_{-2},y_{-1},y_0,y_1,y_2,\dots\}\subset X$, for which the sequence of 
distances $\{\rho(Ty_i,y_{i+1})\}$ for all available indices $i$ satisfies a certain ``smallness'' condition. 
The parts corresponding to non-negative or non-positive indices are referred as forward or 
backward pseudo-trajectories.}

Identifying the indices with moments of time, it is reasonable to think about 
the ``gaps'' $\gamma_{t_i}:=\rho(Ty_{t_i},y_{t_{i+1}})\ne0$ as amplitudes of 
perturbations of the true trajectory at time $t_i$. 
Therefore we introduce the set of ``moments of perturbations'':
$$\cN(\v{y}):=\{t_i:~ \rho(Ty_{t_i}, y_{t_{i+1}})>0, ~i\in \IZ\}$$ 
ordered with respect to their values, i.e. $t_i<t_{i+1}~\forall i$.

Specifying the ``smallness'' condition, a large number of various types of pseudo-trajectories 
are already considered in the literature (see, for example, discussion in \cite{LS,KKO}). 
In this paper we consider only two of them and introduce a new version as well.

\bdef{For a given $\ep>0$ we say that a pseudo-trajectory $\v{y}$ is of 
\begin{itemize}
\item[(U)] {\em uniform} type, if $\rho(Ty_i,y_{i+1})\le\ep$ for all available indices $i$.
\item[(A)] {\em small on average} type, if $\exists N$ such that
    $\frac1{2n+1}\sum\limits_{i=-n}^n\rho(Ty_i,y_{i+1})\le\ep~~\forall n\ge N$. 
    If the backward pseudo-trajectory is finite 
    $\frac1{n+1}\sum\limits_{i=0}^n\rho(Ty_i,y_{i+1})\le\ep~~\forall n\ge N$.
\item[(R)] {\em rare perturbations} type, if the upper density of the set $\cN(\v{y})$ does not exceed $\ep$. 
    Namely, $\limsup\limits_{n\to\infty}\frac1{2n+1} \#(\cN(\v{y})\cap [-n,n]) \le\ep$. 
    Again if the backward pseudo-trajectory is finite, the condition is  
    $\limsup\limits_{n\to\infty}\frac1{n+1} \#(\cN(\v{y})\cap [0,n]) \le\ep$. 
\end{itemize}}

The U-type pseudo-trajectory is the classical one, introduced by G.~Birkhoff \cite{Bi} and 
D.V.~Anosov \cite{An}. The A-type was proposed by M.~Blank \cite{Bl} in order to take 
care about Gaussian type perturbations. The R-type is a new one. It allows to consider 
large but rare perturbations, which are not covered by previous approaches.\footnote{Without 
        the assumption in (A) that the inequality holds for all large enough $n$
        (for true Gaussian perturbations, the probability of this event is zero),  
        the type (R) would belong to (A).} 
A similar notion, called an {\em ergodic pseudo-orbit}\footnote{There are perturbations 
        of various amplitudes here, but the proportion of large ones tends to zero over time. 
        It is worth noting that the name ``ergodic'' seems inappropriate in this context, 
        since it has nothing to do with classical ergodicity.} 
(see, for example, \cite{WOC}), is a mixture of our U and R types.

To simplify notation we will speak about $\ep$-pseudo-trajectories, when the corresponding 
property is satisfied with the accuracy $\ep$. 

\bigskip

The idea of {\em shadowing} in the dynamical systems theory boils down to the question 
is it possible to approximate pseudo-trajectories of a given dynamical system by true trajectories? 
Naturally, the answer depends on the type of approximation.

\bdef{We say that a true trajectory $\v{x}$ {\em shadows} a pseudo-trajectory $\v{y}$ 
with accuracy $\delta$ (notation $\delta$-shadows): 
\begin{itemize}
\item[(U)] {\em uniformly}, if $\rho(x_i,y_i)\le\delta$ for all available indices $i$.
\item[(A)] {\em on average}, if 
      $\limsup\limits_{n\to\infty} \frac1{2n+1}\sum\limits_{i=-n}^n\rho(x_i,y_i)\le\delta$. 
      Again if the backward pseudo-trajectory is finite, the condition is  
      $\limsup\limits_{n\to\infty} \frac1{n+1}\sum\limits_{i=0}^n\rho(x_i,y_i)\le\delta$.
\item[(L)] {\em in the limit}, if $\limsup\limits_{|n|\to\infty}\rho(x_n,y_n)\le\delta$. 
      If the backward pseudo-trajectory is finite, it is enough to have 
      $\limsup\limits_{n\to\infty}\rho(x_n,y_n)\le\delta$. 
\end{itemize}}

The U-type shadowing was originally proposed by D.V.~Anosov \cite{An}, 
while the A-type was introduced\footnote{since pseudo-trajectories with large 
       perturbations cannot be uniformly shadowed.} by M.~Blank \cite{Bl}.
Naturally, the types of pseudo-trajectories and the types of shadowing may be paired 
in an arbitrary way.

\bdef{We say that a DS $(T,X,\rho)$ satisfies the {\em $(\alpha+\beta)$-shadowing property} 
with $\alpha\in\{U,A,R\},~ \beta\in\{U,A,L\}$ if $\forall\delta>0~\exists\ep>0$ such that 
each $\ep$-pseudo-trajectory of $\alpha$-type can be shadowed in the $\beta$ sense with the 
corresponding accuracy $\delta$.}

For example, $(U+U)$ stands for the classical situation of the uniform shadowing of uniformly 
perturbed systems, while $(R+A)$ corresponds to the average shadowing in the case 
of rare perturbations.

There is a large number of open questions related to the shadowing problem and let 
me start with one of them. Under what conditions on the map does the presence of 
a certain shadowing type for each pseudo-trajectory under a single perturbation 
imply one or another type of shadowing property for the system? 
The reason for this question is that the case of a single perturbation is much simpler, 
and therefore the idea of obtaining information about other types of perturbations 
from this fact is quite attractive.
The answer is known (although very partially, see \cite{Ach}) only in the case of 
U-shadowing of the so-called positively expansive\footnote{Roughly speaking,  
    this means that if two forward trajectories are uniformly close enough to each 
    other, then they coincide. In particular, this property is satisfied for expanding maps.} 
dynamical systems, if additionally one assumes that the single perturbation does 
not exceed $0<\ep\ll1$. 

In order to give the answer to this question we introduce the following property.

\bdef{We say that a trajectory $\v{z}$ {\em glues}  together semi-trajectories $\v{x}, \v{y}$ 
with accuracy rate $\phi:\IZ\to\IR_+$ if 
\beq{e:glu}{
     \rho(x_k, z_k)\le\phi(k)\rho(x_0,y_0) ~~\forall k<0, \quad  
     \rho(y_k, z_k)\le\phi(k)\rho(x_0,y_0) ~~\forall k\ge0 .}
In other words $\v{z}$ approximates both the backward part of $\v{x}$ and the 
forward part of $\v{y}$ with accuracy controlled by the rate function $\phi$ 
and the gap between $\v{x}, \v{y}$ at time $t=0$. }

Without loss of generality, we assume that the functions $\phi(|k|)$ and $\phi(-|k|)$ 
are monotonic. Indeed, replacing a general $\phi$ by its monotone envelope 
$$\t\phi(k):=\function{\sup_{i\le k}\phi(i) &\mbox{if } k<0 \\  
                               \sup_{i\ge k}\phi(i) &\mbox{if } k\ge0} ,$$ 
we get the result. 

\bdef{We say that the DS $(T,X,\rho)$ satisfies the {\em gluing property} with the rate-function 
$\phi:\IZ\to\IR$ (notation $T\in G(\phi)$) if for any pair of trajectories $\v{x}, \v{y}$ there is 
a trajectory $\v{z}$, which glues them at time $t=0$ with accuracy $\phi$.}

\begin{remark}
If $T\in G(\phi)$, then $\forall \tau\in\IZ$ for any pair of trajectories $\v{x}, \v{y}$ 
there exists a trajectory $\v{z}$, which glues them at time $t=\tau$ with accuracy $\phi$. 
\end{remark}
Indeed, for a given $\tau$ consider a pair of trajectories $\v{x'}, \v{y'}$ obtained from 
$\v{x}, \v{y}$ by the time shift by $\tau$, namely $x'_i:=x_{i+\tau}, ~y'_i:=y_{i+\tau}, ~\forall i$. 
Then since $\v{x'}, \v{y'}$ may be glued together at time $t=0$ with accuracy $\phi$,  
we deduce the same property for $\v{x}, \v{y}$ by the time $t=\tau$. \qed

This property may be explained in terms similar to those which are actively used in the theory 
of smooth hyperbolic dynamical systems. Denote by $\v{x}^-:=\{\dots,x_{-2},x_{-1},x_0=x\}$ and 
$\v{x}^+:=\{x=x_0,x_1,x_x,\dots\}$ backward and forward semi-trajectories of the point $x\in X$, 
and consider the sets:
$$ W^-(\v{x}^-):=\{\v{z}\subset X:~~\rho(x_k,z_k)\toas{k\to\-\infty}0\} .$$ 
$$ W^+(\v{x}^+):=\{\v{z}\subset X:~~\rho(x_k,z_k)\toas{k\to\infty}0\} .$$ 
Then the gluing property means that for each pair of semi-trajectories $\v{x}^-$ and $\v{y}^+$ 
the sets $W^+(\v{x}^+)$ and $W^-(\v{y}^-)$ have a non-empty intersection. Additionally 
the rate of convergence in the definition of the sets $W^\pm$ is controlled by the rate 
function $\phi$.

The origin of the gluing property is the so-called local product structure, introduced by D.V.~Anosov 
for uniformly hyperbolic DS. The local product structure means that for a pair 
of close enough points their stable and unstable manifolds intersect, and the orbit of 
the point of intersection approximates the corresponding semi-trajectories with 
an error exponentially decreasing in time. In our notation this means $\phi(k)=Ce^{-b|k|}$. 
Later in \cite{Bl} this property has been extended to the global one, but only for 
the uniformly hyperbolic DS.

\begin{remark} (Necessity) The gluing property is necessary for the A and R types of shadowing, 
but not for the L type.
\end{remark}

Already in the simplest case of a single large perturbation, the average shadowing 
implies the gluing of any pair of semi-trajectories. If $D:={\rm diam}(X)<\infty$ we may choose 
$\phi\equiv D$, otherwise, assuming that the perturbations are bounded by a constant 
$D<\infty$, one shows (see Section~\ref{s:main}) that during the gluing procedure 
the gaps between the glued segments of true trajectories cannot exceed $KD$, 
where $K$ depends only on $T$, but not on the particular pseudo-trajectory. 
Therefore in the unbounded case we set $\phi\equiv KD$. Contrary to this in the case of 
limit shadowing (type L) we cannot control deviations between the pseudo-trajectory 
and the true trajectory. \qed

\begin{remark} 
In order to obtain meaningful applications of the gluing property, additional assumptions 
about the rate-function $\phi$ must be made. In what follows, we will only assume  
summability of this function: 
\beq{e:sum}{ \Phi:=\sum_k\phi(k)<\infty .}
\end{remark}

In distinction to the gluing property itself, the summability does not follow from 
either the uniform shadowing or from from the average one. 
Indeed, both these types of shadowing imply only that 
$$ \limsup_{n\to\infty} \frac1{2n+1}\sum_{k=-n}^n \phi(k) \le \delta ,$$
which is insufficient for the summability. 

Moreover, even if the rate-function $\phi$ is summable, its monotone envelope 
$\t\phi$ might be non-summable. Consider an example:
$$ \phi:~~~ \cdots k^{-2}~\underbrace{0\dots0}_{k-1}\cdots 3^{-2}~0~0~ 2^{-2}~0~1~ 
                    0~2^{-2} ~0~0~3^{-2} \cdots k^{-2}\underbrace{0\dots0}_{k-1} \cdots $$
Here $\phi$ is an even function, vanishing at $\pm\infty$ with $\phi(0)=1$. Clearly, 
$\phi$ is summable: $\sum_i\phi(i)=\pi^2/3<\infty$. On the other hand, $\phi(|k|)$ 
is non-monotonic and its monotone envelope $\t\phi$ is no longer summable: 
$\sum_i\t\phi(i) = 1+ 2\sum_{k>1} k\cdot k^{-2} = 1+ 2\sum_{k>1} k^{-1} = \infty$.

\section{Main result}\label{s:main}
\begin{theorem}\label{t:main}
Let $T:X \to X$ be a map from a metric space $(X,\rho)$ into itself, and 
let $T\in G(\phi)$ with $\Phi:=\sum_k\phi(k)<\infty$. Then
\begin{itemize}
\item[(a)] $T\in (U+U)$,
\bwait2{
\item[(b)] $T\in (R+A)$ if additionally the perturbations are bounded by some constant $D<\infty$.}
\end{itemize}
\end{theorem}

\begin{remark}\label{r:bounded}
\begin{enumerate}
\item in the case (U+U) it is enough to check the gluing property for $\v{x}, \v{y}$ with 
        $\rho(x_0,y_0)\le\ep_0\ll1$,
\item in the case (R+A) it is enough to check the gluing property for $\v{x}, \v{y}$ with 
        $\rho(x_0,y_0)\le De^\Phi$,
\item if $\diam(X):=\sup_{u,v\in X}\rho(u,v)<\infty$, then the perturbations cannot exceed 
        $\diam(X)$ and in the case (R+A) the gluing property may be weaken to 
\beq{e:glu-w}{\rho(x_k, z_k)\le\phi(k) ~~\forall k<0, \quad  \rho(y_k, z_k)\le\phi(k) ~~\forall k\ge0,} 
       i.e. the factor $\rho(x_0,y_0)$ is dropped compared to (\ref{e:glu}).
\end{enumerate} 
\end{remark}

\begin{remark}
We expect more sophisticated estimates (compared to the present proof) 
to give $T\in (A+A)$ under the same assumptions. 
\end{remark}

\proof In fact, in each case, we will prove a stronger ``linear'' version of shadowing, 
namely that there is a constant $K=K(\phi)<\infty$, such that for each $\ep>0$ small enough 
for each $\ep$-pseudo-trajectory there is a true trajectory approximating it 
(in the corresponding sense) with accuracy $\delta\le K\ep$ . 

The proof goes as follows. For an $\ep$-pseudo-trajectory under rare perturbations 
$\v{y}:=\{y_i\}$ consider in detail the set $\cN(\v{y})$ of moments of perturbations 
$\dots < t_{-2} < t_{-1} < t_0 < t_1 < t_2<\dots$. 
Between the moments of time $t_k$ there are no perturbations and hence $\v{y}$ can 
be divided into segments of true trajectories. Thanks to the $G(\phi)$ property each pair of 
consecutive segments of true trajectories can be ``glued'' together by a true trajectory with the 
controlled accuracy. 

\begin{figure}\begin{center}
\begin{tikzpicture}[scale=0.75]
      \draw [-](-.8,0) to (16.5,0);
      %\draw (-.5,0) to (-.5,-.2); \node at (-.5,-.5){$t_{-9}$};
      \draw (.5,0) to (0.5,-.2); \node at (.5,-.5){$t_{-8}$};
      \draw (1.5,0) to (1.5,-.2); \node at (1.5,-.5){$t_{-7}$};
      \draw (2.5,0) to (2.5,-.2); \node at (2.5,-.5){$t_{-6}$};
      \draw (3.5,0) to (3.5,-.2); \node at (3.5,-.5){$t_{-5}$};
      \draw (4.5,0) to (4.5,-.2); \node at (4.5,-.5){$t_{-4}$};
      \draw (5.5,0) to (5.5,-.2); \node at (5.5,-.5){$t_{-3}$};
      \draw (6.5,0) to (6.5,-.2); \node at (6.5,-.5){$t_{-1}$};
      \draw (7.2,0) to (7.2,-.2); \node at (7.2,-.5){$t_{0}$};
      \draw (8.,0) to (8.0,-.2); \node at (8.0,-.5){$t_{1}$};
      \draw (9.,0) to (9.0,-.2); \node at (9.0,-.5){$t_{2}$};
      \draw (10.,0) to (10.0,-.2); \node at (10.0,-.5){$t_{3}$};
      \draw (11.,0) to (11.0,-.2); \node at (11.0,-.5){$t_{4}$};
      \draw (12.,0) to (12.0,-.2); \node at (12.0,-.5){$t_{5}$};
      \draw (13.2,0) to (13.2,-.2); \node at (13.2,-.5){$t_{6}$};
      \draw (14.2,0) to (14.2,-.2); \node at (14.2,-.5){$t_{7}$};
      \draw (15.2,0) to (15.2,-.2); \node at (15.2,-.5){$t_{8}$};
      \draw (16.2,0) to (16.2,-.2); \node at (16.2,-.5){$t_{9}$};
      \node at (-.5,.5){$n=1$}; \node at (-.5,1.0){$n=2$}; \node at (-.5,1.5){$n=3$};
      %\draw [thick,<->] (6.5,.5) to (8.0,.5); %\node at (1.5,-.5){$t_{-1}$};
      \draw [line width=1.5pt,,<->] (6.5,.5) to (8.0,.5); 
      \draw [line width=1.5pt,,<->] (8.0,.5) to (10.0,.5);
      \draw [line width=1.5pt,,<->] (10.0,.5) to (12.0,.5);
      \draw [line width=1.5pt,,<->] (12.0,.5) to (14.2,.5);
      \draw [line width=1.5pt,,<->] (14.2,.5) to (16.2,.5);
      \draw [line width=1.5pt,,<->] (4.5,.5) to (6.5,.5); 
      \draw [line width=1.5pt,,<->] (2.5,.5) to (4.5,.5); 
      \draw [line width=1.5pt,,<->] (.5,.5) to (2.5,.5); 
      \draw [line width=1.5pt,,<->] (6.5,1.0) to (10.0,1.0); 
      \draw [line width=1.5pt,,<->] (10.0,1.0) to (14.2,1.0); 
      \draw [line width=1.5pt,,<->] (2.5,1.0) to (6.5,1.0); 
      \draw [line width=1.5pt,,->] (.5,1.0) to (2.5,1.0); 
      \draw [line width=1.5pt,,<-] (14.2,1.0) to (16.2,1.0); 
      \draw [line width=1.5pt,,<->] (2.5,1.5) to (10.0,1.5); 
      \draw [line width=1.5pt,,->] (.5,1.5) to (2.5,1.5); 
      \draw [line width=1.5pt,,<-] (10.0,1.5) to (16.2,1.5); 
\end{tikzpicture}\end{center}
\caption{Order of the parallel gluing.}\label{f:sgluing} 
\end{figure}
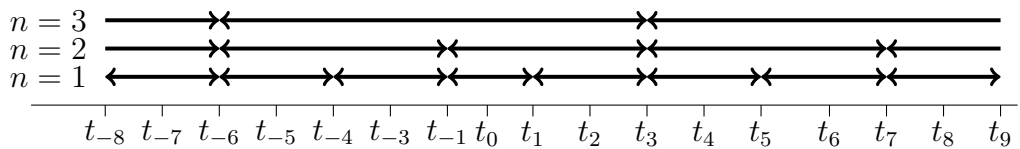

In our construction (see Fig.~\ref{f:sgluing}) we first simultaneously ``glue'' together pairs 
of segments around the moments of perturbations $t_i$ with even indices: $i_{\pm2k}$, 
getting longer segments of the ``gluing'' true trajectories. 
On each next step we simultaneously ``glue'' together consecutive pairs of already 
obtained segments. Therefore, on the $n$-th step of the construction we obtain a new 
pseudo-trajectory $\v{z}^{(n)}$, consisting of half (compared to $\v{z}^{(n-1)}$) the number 
of segments of true trajectories (i.e. only a half of the moments of perturbations remain) 
with an exponentially growing lengths, but with possibly larger gaps between them 
(compared to original gaps). In the limit we get the approximation of the entire pseudo-trajectory. 

The procedure above can be called a parallel gluing construction. 
Another possibility is to use a consecutive gluing, which can be described 
as follows. Starting from some segment of the trajectory (say, between 
the time moments from $t_0$ to $t_1$), we glue it first with the right 
neighbor, then with the left one (or vice versa). The construction of 
sequential gluing consists in sequential gluing each time of a new 
segment of the trajectory to the already glued ones.
In fact, the construction used in \cite{Bl} to prove the average shadowing property 
for Anosov systems in the above terminology is exactly the consecutive gluing. 
The advantage of the consecutive gluing construction is that the corresponding 
calculations are much simpler, but on closer examination it turns out that 
in order to apply them, it is necessary to make much stronger assumptions about  
the rate function $\phi$, in particular, that $\phi(\pm1)<1$. 
Even for uniformly hyperbolic systems, this can be done only for the so-called 
Lyapunov metric $\rho$, and not for the general one. 
In most of the examples discussed in Section~\ref{s:ver}   
the value $\phi(\pm1)$ turns out to be quite large. 

To estimate the approximation errors we find the accuracy of the gluing of a pair of 
segments of true trajectories: $v_{-N^-}, v_{-N^-+1},\dots,v_{-1}$ and $v_0,v_1,\dots,v_{N^+}$. 
By the $G(\phi)$-property there is a trajectory $\v{z}\subset X$ such that 

$$ \rho(v_k, z_k)\le\phi(k)\rho(Tv_{-1},v_0) ~~\forall k\in \{-N^-,\dots,N^+\}$$
Therefore 
$$ \sum_{k=-N^-}^{N^+}\rho(z_k, v_k) \le \rho(Tv_{-1},v_0) \sum_{k}\phi(k) 
    = \Phi\cdot \rho(Tv_{-1},v_0) .$$
There are three important points here: 
\begin{enumerate}
\item the result depends only on the gap $\rho(Tv_{-1},v_0)$ between the ``end-points'' 
        of the glued segments of trajectories;
\item after the gluing of a pair of segments of trajectories the gaps between the end-points 
        on the next step of the procedure may become larger than the original gaps in $\v{y}$;
\item each moment of perturbation $t_i$ is taken into consideration only once in the process of gluing.
\end{enumerate}

Let us first prove part (b) of the assertion.

\begin{figure}\begin{center}
\def\G#1#2#3{\draw [thick] (1.6+#1,0+#2) to (1.6+#1,1.5+#2); \node at (1.6+#1,-.5+#2){#3};
      \draw [thick] (0+#1,0+#2) .. controls (1+#1,.5+#2) and (1.5+#1,2.5+#2) .. (2+#1,1+#2) 
                                 .. controls (2.1+#1,.5+#2) and (2.5+#1,.3+#2) .. (3+#1,0+#2);}
\begin{tikzpicture}[scale=0.75]
      \draw [-](0,0) to (11,0);  \node at (5.1,.8){$D$};
      \G{1}{.05}{$t_{-1}$}; \G{2}{.05}{$t_{0}$}; \G{3.8}{.05}{$t_{1}$};
      \G{6}{.05}{$t_{2}$}; \G{7.5}{.05}{$t_{3}$};
\end{tikzpicture}\end{center}
\caption{Contributions to the upper bound of the gluing error.}\label{f:error} 
\end{figure}
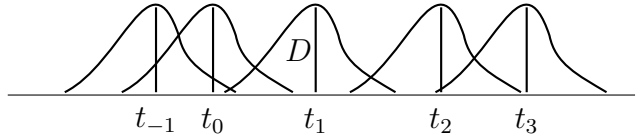

If $D:=\diam(X)<\infty$ the gaps cannot exceed $D$. Therefore the contributions from  
each gap to the final approximation error are simply adding together (see Fig.~\ref{f:error}).

Using this, we get
$$ \rho(z_t^{(n)},y_t)  \le D \sum_i \phi(t-t_i) \le D \Phi \quad\forall t\in\IZ.$$

Similarly, for each $k>0$ 
$$ \rho(z_t^{(n)},z_t^{(n+k)})  \le D \sum_{|j|>r(t,n)} \phi(j) \toas{n\to\infty}0 ,$$
where $r(t,n)\toas{n\to\infty}\infty$. Therefore for any given $t$ the sequence 
$\{z_t^{(n)}\}$ is fundamental and converges as $n \to\infty$ to the limit $z_t$, 
where $\v{z}:=\{z_t\}$ is the true trajectory of our system.

The averaged approximation error can be estimated as follows:
$$ Q_m(\v{z}^{(n)},\v{y}):=\frac1{2m+1}\sum_{|k|\le m}\rho(z_k^{(n)},y_k) 
    \le \Phi D \cdot \frac{\#(\cN(\v{y})\cap [-m,m])}{2m+1} .$$
Hence 
$$ \limsup_{m\to\infty} Q_m(\v{z}^{(n)},\v{y}) \le \Phi D\ep ,$$
which implies the same inequality for $\v{z}$ and hence the average shadowing 
in the case under question. 

Note that in these calculations we always use only the largest possible amplitudes 
of perturbations $D$. This effectively reduces the gluing condition (\ref{e:glu}) to 
its weak version (\ref{e:glu-w}) (without the dependence on the actual value of 
the perturbation), formulated in Remark~\ref{r:bounded}.

\bigskip

To consider the unbounded case $\diam(X)=\infty$  we need the following simple 
statement about infinite products. 

\begin{lemma}\label{e:prod} For any sequence of real numbers $\{b_k\}_{k\ge1}$ we have 
$$ \limsup_{n\to\infty}\prod_{k\ge1}^n(1+b_k) \le e^{\limsup\limits_{n\to\infty}\sum_{k=1}^n b_k} .$$
If additionally $b_k\ge0~~\forall k\in\IZ_+$, then
$$ \liminf_{n\to\infty}\prod_{k\ge1}^n(1+b_k) \ge 1 + \liminf_{n\to\infty}\sum_{k=1}^n b_k .$$
\end{lemma}
\proof Denote $Q_n:=\prod_{k\ge1}^n(1+b_k), ~S_n:=\sum_{k=1}^n b_k$. 
We proceed by induction on $n$. To check the basic step with $n=1$, let us show that 
$$1+v \le e^v~\forall v\in\IR .$$ 
For $v=0$ we have equality, while 
$$ \frac{d}{dv} (1+v) = 1 <  e^v = \frac{d}{dv}e^v ~~\forall v>0,$$
$$ \frac{d}{dv} (1+v) = 1 >  e^v = \frac{d}{dv}e^v ~~\forall v<0.$$
Therefore the graph of $e^v$ lies above the straight line $1+v$ with the only 
tangent point at $v=0$. 

This implies that $Q_1\le e^{S_1}$. Assume that  $Q_n\le e^{S_n}$ for some $n\in\IZ_+$ 
and prove the same inequality for $n+1$. We get 
$$ Q_{n+1}=(1+b_{n+1})Q_n \le (1+b_{n+1}) e^{S_n} \le e^{b_{n+1}}e^{S_n} = e^{S_{n+1}} .$$
Similarly but even simpler we get the estimate from bellow if $b_k\ge0~\forall k$. 
The basic step here is trivial, while the induction step can be done as follows:
$$ Q_{n+1}=(1+b_{n+1})Q_n \ge (1+b_{n+1})(1+S_n) = 1 + S_n + b_{n+1} + b_{n+1}S_n 
    \ge 1 + S_{n+1} .$$
Passing to the upper/lower limits as $n\to\infty$ we get the result. \qed

Now we are ready to proceed with the unbounded case $\diam(X)=\infty$. 
Here we use the additional assumption about the perturbations amplitudes  $D:=\sup_k\rho(Ty_k,y_{k+1})<\infty$ 
to take care about the growth of the gaps during the gluing. 

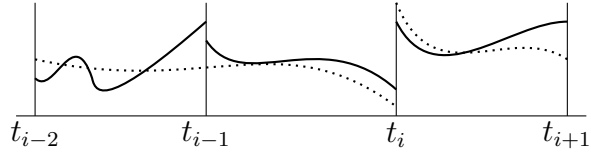
\begin{figure}\begin{center}
\begin{tikzpicture}[scale=0.5]
      \draw [-](0,0) to (15,0);
      \draw[thin] (.5,0) to (.5,3); \node at (.5,-.5){$t_{i-2}$};
      \draw[thin] (5,0) to (5,3); \node at (5,-.5){$t_{i-1}$};
      \draw[thin] (10,0) to (10,3); \node at (10,-.5){$t_{i}$};
      \draw[thin] (14.5,0) to (14.5,3); \node at (14.5,-.5){$t_{i+1}$};
      \draw [thick] (.5,1) .. controls (1,.5) and (1.5,2.5) .. (2,1) 
                                 .. controls (2.1,.5) and (2.5,.3) .. (5,2.5);  %Bezier
      \draw [thick] (5,2) .. controls (6,.5) and (8,2.5) .. (10,.7);  %Bezier 
      \draw [thick] (10,2.5) .. controls (11,.5) and (13,2.5) .. (14.5,2.5);  %Bezier 
      \draw [thick,dotted, black] (.5,1.5) .. controls (5,.5) and (7,2.5) .. (10,.25);  %Bezier 
      \draw [thick,dotted, black] (10,3) .. controls (11,.5) and (13,2.5) .. (14.5,1.5);  %Bezier 
      %\put(45,1){\draw [thick] (1.5,0) to (1.5,1.5); \node at (1.5,-.5){$t_{-1}$};
      %\draw [thick] (0,0) .. controls (1,.5) and (1.5,2.5) .. (2,1) 
       %                          .. controls (2.1,.5) and (2.5,.3) .. (3,0);}  %Bezier
%      \draw [thick]  (1,10) .. controls (3,7.8) and (3.5,6) .. (3,5.2);  %Bezier
\end{tikzpicture}\end{center}
\caption{Gluing a pair of segments of true trajectories around $t_{i-1}$. The glued 
segments are indicated by dotted lines.}\label{f:gluing} 
\end{figure}

On the $n$-th step of the gluing process we have a pseudo-trajectory $\v{z}^{(n)}$ with 
a bi-infinite collection of gaps $\{\gamma_{t_i}^{(n)}\}_{i\in\IZ}$, occurring at certain 
moments of time $t_i$ (see Fig.~\ref{f:gluing}). 

By the gluing property, when gluing segments of true trajectories around the point 
$t_{i-1}$ at the $n$-th step of the procedure, we obtain a recursive upper bound 
for the gap at point $t_i$ on $(n+1)$-th step:
\beq{e:rec-est0}{
  \gamma_i^{(n+1)} \le \gamma_i^{(n)} 
                                   + \phi_-^{(n)}\gamma_{i-1}^{(n)}  + \phi_+^{(n)} \gamma_{i+1}^{(n)} ,}
where $\phi_-^{(n)}=\phi(t_i-t_{i-1})$ and $\phi_+^{(n)}=\phi(t_i-t_{i+1})$. 
Indeed, $\phi_-^{(n)}\gamma_{i-1}^{(n)}$ is the upper bound for the approximation error coming 
from the left, while $\phi_+^{(n)} \gamma_{i+1}^{(n)}$ is the upper bound for the 
approximation error coming from the right. 

Denote by $\gamma^{(n)}:=\sup_i \gamma_{t_i}^{(n)}$ the maximal value of the gaps on the $n$-the 
step of the procedure, and by $\tau^{(n)}:=\inf_i|t_i^{(n)} - t_{i+1}^{(n)}|$ -- the length 
of the smallest segment of true trajectories. Then using (\ref{e:rec-est0}) and the monotonicity 
of the functions $\phi(\pm|k|)$ we get 
$$ \gamma^{(n+1)} \le \gamma^{(n)} 
                   + \phi(\tau^{(n)})\gamma^{(n)}  + \phi(-\tau^{(n)}) \gamma^{(n)} 
     = \gamma^{(n)} \cdot \left(1 + \phi(-\tau^{(n)}) + \phi(\tau^{(n)}) \right) .$$
Continuing this and passing from $n$ to $n-1$, etc., we obtain the following estimate:
$$  \gamma^{(n+1)} \le \gamma^{(0)} \cdot 
           \prod_{k=0}^n \left(1 + \phi(-\tau^{(k)}) + \phi(\tau^{(k)}) \right) .$$

By Lemma~\ref{e:prod} the right hand side in the last inequality may be estimated 
from above  
$$  \gamma^{(n+1)} \le  \gamma^{(0)} \cdot 
                 \exp\left(\sum_{k\ge0} \left(\phi(-\tau^{(k)}) + \phi(\tau^{(k)}) \right) \right) 
                \le D \exp\left(\sum_k \phi(k) \right) = De^\Phi ,$$
since the function $\tau^{(k)})$ is strictly increasing, and in the case under consideration 
$\gamma^{(0)} \le D$. 

Thus the gaps $\gamma_{t_i}^{(n)}$ are uniformly bounded from above by $De^\Phi$. 
Applying the same trick, as in the bounded case (but using $De^\Phi$ instead of $D$), 
we get a slightly worse estimate which is still sufficient for our purpose:
$$ Q_m(\v{z}^{(n)},\v{y}):=\frac1{2m+1}\sum_{|k|\le m}\rho(z_k^{(n)},y_k) 
 \le De^{\Phi}\cdot \frac{\#(\cN(\v{y})\cap [-m,m])}{2m+1} .$$
Hence 
$$ \limsup_{m\to\infty} Q_m(\v{z}^{(n)},\v{y}) \le De^{\Phi}\ep .$$

Finally, the same reasoning as in the first part of the proof is used to verify 
the convergence of pseudo-trajectories $\v{z}^{(n)}$ to the true trajectory $\v{z}$, 
and to check the inequality above for $\v{z}$ instead of $\v{z}^{(n)}$.

\bigskip

It remains to prove part (a) of the assertion: $T\in (U+U)$. Here, despite the perturbations 
are uniformly small, they can occur arbitrary often. Therefore, one cannot make use 
of large distances between them, as in the cases already considered. 
Nevertheless, using the same argument as in the second part of the proof 
(replacing the constant $D$ by $\ep$), 
one gets the upper bound for the maximal possible value of gaps between the 
segments of trajectories in the process of gluing. Now instead of the upper bound for 
the average in time approximation error one needs to get the upper bound for the 
uniform approximation. 

Since the contribution of each gap to the final uniform approximation error 
is summed (see Fig.~\ref{f:error}), there is a finite number $K\ne K(\ep)$ such that
$$ \rho(z_t,y_t)  \le \ep K \sum_i \phi(t-t_i) \le \ep K \Phi \quad\forall t\in\IZ.$$
The proof of the theorem is complete. \qed 

\begin{remark}
The global version of the gluing property is used only for the 1-st part 
of the proof. For the remaining parts, it suffices to consider pairs of trajectories 
$\v{x}, \v{y}$, for which $\rho(x_0, y_0)\le KD$. 
This greatly simplifies the verification of the gluing property in these cases. 
\end{remark}

\begin{remark}
Our estimates are very crude, because we assume only that the 
perturbations are uniformly bounded by $D$. In fact, using information about their 
actual size or the distribution of their lengths one might try to get a better result. 
This is what we had in mind writing about more sophisticated estimates (immediately 
after the formulation of Theorem~\ref{t:main}). 
Unfortunately, this is not possible within the approach used in the present proof because 
during the gluing procedure we cannot keep track of individual gaps. 
\end{remark}

\section{Verification of the gluing property}\label{s:ver}

\begin{example}\label{e:hyp} (Hyperbolic affine mapping) 
Let $X:=\IR^2, Tx:=Ax + a$, where the matrix $A$ has real eigenvalues $\la_1>1>\la_2>0$, 
corresponding to non-collinear unit eigenvectors $e_1,e_2\in\IR^2$, and $a\in\IR^2$. 
\end{example}

\begin{figure}\begin{center}
\begin{tikzpicture}[scale=0.5]
     \draw [-](0,0) to (10,0); \draw [-](2,-1) to (2,8);
     \draw [-](0,-1) to (10,2.5); \draw [<->](0,-1) to (5,.75); \node at (5,.2){$e_1$}; 
     \draw [-](3,-1) to (6,8);     \draw [>-<](3,-1) to (4,2);    \node at (4.4,1.5){$e_2$}; 
     \node at (10,1.9){$L_1$};  \node at (3.6,-.65){$L_2$}; 
     \draw [thick]  (0,-.8) .. controls (2,0) and (3.5,1) .. (5.8,8);  %Bezier
     \draw [thick]  (6.2,8) .. controls (3.5,0) and (4.5,1) .. (10,2.7);  %Bezier
     \draw [fill] (5.5,1.5) circle (0.2cm);  \node at (6.2,2){$y_0$}; 
     \draw [fill] (3.6,2.5) circle (0.2cm);  \node at (2.8,2.5){$x_0$};
     \draw [dotted](3.6,2.5) to (7,4); \draw [dotted](5.5,1.5) to (6.5,4.3);
     \draw [fill] (6.3,3.7) circle (0.2cm); \node at (7,3.7){$z_0$}; 
     \draw [thick, dotted] (6.4,8) .. controls (5,3.3) and (7,2.1) .. (10,2.9);  %Bezier
%      \draw [thick] (0,0) .. controls (1,.5) and (1.5,2.5) .. (2,1) 
%                                 .. controls (2.1,.5) and (2.5,.3) .. (3,0);}  %Bezier
%      \draw [thick]  (1,10) .. controls (3,7.8) and (3.5,6) .. (3,5.2);  %Bezier
\end{tikzpicture}\end{center}
\caption{Hyperbolic affine mapping. Gluing of $\v{x}$ and $\v{y}$. Two typical trajectories are 
             indicated by thick lines, while the gluing trajectory by a thick dotted line.}
\label{f:hyp-ex}  \end{figure}
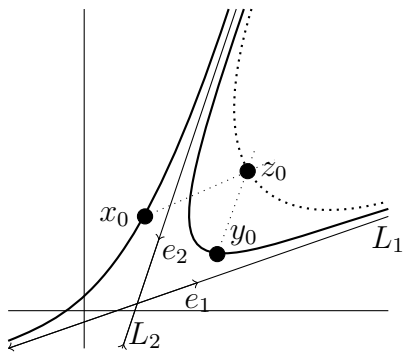

\begin{proposition} \label{p:hyp}
Trajectories $\v{x},\v{y}$ of the Hyperbolic affine mapping $T:\IR^2\to\IR^2$ in the example~\ref{e:hyp}  
satisfies the $G(\phi)$ property with $\phi(k):=\function{C\la_1^{-k} &\mbox{if~} k\ge0 \\
                                                                                 C\la_2^{|k|} &\mbox{if~} k\le0 }$  
if they belong to the same half-plane bounded by one of the lines $L_i$.
The metric $\rho$ in $\IR^2$ is assumed to be induced by a norm. 
\end{proposition}
\proof We need to check that for any backward semi-trajectory $\v{x}:=\{\dots,x_{-1},x_0\}$ 
and a forward semi-trajectory $\v{y}:=\{y_0,y_1,\dots\}$ of the map $T$ there is 
a trajectory $\v{z}:=\{\dots,z_{-1},z_0,z_1,\dots\}$ such that 
$\rho(z_n,x_n)\le C\la_2^{|n|}~~\forall n\le0$ and $\rho(z_n,y_n)\le C\la_1^{-n}~~\forall n\ge0$.

The $n$-th point ($n\in\IZ$) of a trajectory $\v{v}$ under the action of $T$ is uniquely written as 
$$ v_n=\la_1^n\alpha(\v{v})e_1 + \la_2^n\beta(\v{v})e_2 + A^na ,$$
where $v_0=(\alpha(\v{v})e_1, \beta(\v{z})e_2)$. 

Therefore by the triangle inequality for any $\v{z}$ 
$$ \rho(z_n,x_n) \le \la_1^n|\alpha(\v{z}) - \alpha(\v{x})| 
                           + \la_2^n|\beta(\v{z}) - \beta(\v{x})| ~~\forall n\le 0 ,$$
$$ \rho(z_n,y_n) \le \la_1^n|\alpha(\v{z}) - \alpha(\v{y})| 
                           + \la_2^n|\beta(\v{z}) - \beta(\v{y})| ~~\forall n\ge 0 .$$
Similarly, 
$$ \rho(z_n,x_n) \ge -\la_1^n|\alpha(\v{z}) - \alpha(\v{x})| 
                           + \la_2^n|\beta(\v{z}) - \beta(\v{x})| ~~\forall n\le 0 ,$$
$$ \rho(z_n,y_n) \ge \la_1^n|\alpha(\v{z}) - \alpha(\v{y})| 
                           - \la_2^n|\beta(\v{z}) - \beta(\v{y})| ~~\forall n\ge 0 .$$
Therefore, the distances $\rho(z_n,x_n)$ and $\rho(z_n,y_n)$ are 
uniformly bounded in $n$ if and only if  
$$ \alpha(\v{z}) = \alpha(\v{x}), ~~ \beta(\v{z}) = \beta(\v{x}) .$$
This is achieved as follows. Consider two straight lines $x_0+s_1e_1,~s_1\in\IR$ 
and $y_0+s_2e_2,~s_1\in\IR$ (see Fig.~\ref{f:hyp-ex}). Since $e_i$ are non-collinear, 
there is the only one intersection point $z_0$. Using again non-collinearity, we deduce that 
$$ \max\{|\alpha(\v{z}) - \alpha(\v{x})|, ~  |\alpha(\v{z}) - \alpha(\v{y})| \} \le C\rho(x_0,y_0) ,$$
where the constant $C$ depends only on the angle between the eigenvectors $e_i$. 

Thus the trajectory $\v{z}$ with $(\v{z})_0=z_0$ satisfies the properties above. \qed

\begin{remark}
The important point is that the constant $C=C(A,\rho)$ can be arbitrary large. 
Generalizations to the case $\IR^d,~d>2$ are straightforward, but here the constant $C$ 
depends additionally on the dimension $d$ and grows to infinity with it. 
\end{remark}

\begin{example}\label{e:Anosov} (Anosov diffeomorphism)
Let $X:=\Tor^2$ be a unit 2-dimensional torus and let $T:X\to X$ be a uniformly hyperbolic
diffeomorphism. 
\end{example}

The simplest map that satisfies the above properties is $Tx:=Ax \mod1$, where $A$ is an 
integer matrix with the determinant equal 1 on modulus. For exact definition of the 
uniformly hyperbolic system we refer the reader to numerous publications to the subject 
(see, for example, \cite{An2, Bo, BPSJ, Bl24}).

\begin{proposition} \label{p:Anosov} 
For the map $T:\Tor^2\to\Tor^2$ of the example~\ref{e:Anosov} 
there exists a special (Lyapunov) metric $\rho$ and $\la>1$, for which this system 
satisfies the $G(\phi)$ property with $\phi(k):=e^{-\la|n|}~~\forall n\in\IZ$.   
\end{proposition}

This result follows from the global product structure for a hyperbolic system proven in \cite{Bl}. 
The local version of this property, which asserts the intersection of stable and unstable 
local manifolds of sufficiently close points, is well known (see, for example,  \cite{An2, Bo, BPSJ}). 
In the global version, the locality assumption is dropped. 

It is worth noting that none of these results follow from Proposition~\ref{p:hyp}, 
and not vice versa. Indeed, the diffeomorphism under study is nonlinear and 
the proof of Proposition~\ref{p:Anosov} is based on the construction of 
arbitrary thin Markov partitions and their mixing properties. On the other 
hand, these constructions fail in the example~\ref{e:hyp}. Moreover, 
Proposition~\ref{p:hyp} holds for an arbitrary metric $\rho$ (induced by a norm), 
while Proposition~\ref{p:Anosov} holds only for a special (Lyapunov) metric.

\bigskip

Now we are ready to turn to discontinuous and non invertible mappings. 

\bdef{A map $T:X\to X$ is said to be {\em piecewise bijective (PB)} if there is 
a partition $\{X_i\}$ of $X$ such that $T:X_i\to TX_i$ is a bijection $\forall i$. 
Each point $u\in X$ belongs to a single element of the partition, which we 
denote by $X_{i(u)}$. For a given point $v\in V$ the inverse map 
$T_v^{-1}:=T_{|X_{i(v)}}^{-1}$ with respect to $v$ is well defined at $TX_{i(v)}$.}

The maps in all examples bellow are of type PB. 

\begin{example}\label{e:plin} (Piecewise linear mapping)
$X:=[0,1]$, $Tx:=\function{ax &\mbox{if } x<c \\  bx + 1-b  &\mbox{otherwise}}$. 
\end{example}

%\bcr{Another version of ``right'' map: $\max(a,b)<1$ -- check is it enough?}

\begin{figure} \begin{center}
\begin{tikzpicture}[scale=0.75]
\put(20,0){
     \draw [-](0,0) to (10,0) to (10,10) to (0,10) to (0,0); 
     \draw [dotted] (0,0) to (10,10); \draw [dotted] (6,0) to (6,10); \node at (6,-0.5){$c$};
     \draw (0,9) to (-0.2,9);\node at (-.5,9){$ac$}; 
     \draw (10,1.5) to (10.2,1.5);\node at (12+.5,1.5){$(1-b)+bc$}; 
     \node at (0.2,-0.5){$0$}; \node at (10,-0.5){1}; 
     \draw [line width=1pt,->](0,0) to (6,9); \draw [line width=1pt,-](10,10) to (6,1.5); 
       \bcb{       
         \draw[fill] (5.5,5.5) circle(.1); \node at (5.5,6){$x_0$};
         \draw [dotted] (5.5,5.5) to (3.7,5.5) to (3.7,3.7);
         \draw[fill] (3.7,3.7) circle(.1); \node at (3.7,4.2){$x_{-1}$};
         \draw [dotted] (3.7,3.7) to (2.5,3.7) to (2.5,2.5);
         \draw[fill] (2.5,2.5) circle(.1); \node at (2.5,3){$x_{-2}$}; }
       \bcb{       % \bcr{  
         \draw[fill] (6.2,6.2) circle(.1); \node at (6.2,6.7){$y_0$};
         \draw [dotted] (6.2,6.2) to (6.2,2.1) to (2.1,2.1);
         \draw[fill] (2.1,2.1) circle(.1); \node at (2.1,1.6){$y_1$};
         \draw [dotted] (2.1,2.1) to (2.1,3.2) to (3.2,3.2);
         \draw[fill] (3.2,3.2) circle(.1); \node at (3.2,2.7){$y_2$};
             } }
\end{tikzpicture}\end{center}
\caption{A discontinuous map, satisfying the gluing property if $ac=b(1-c)=1$.}\label{f:gluing-ex}
\end{figure}
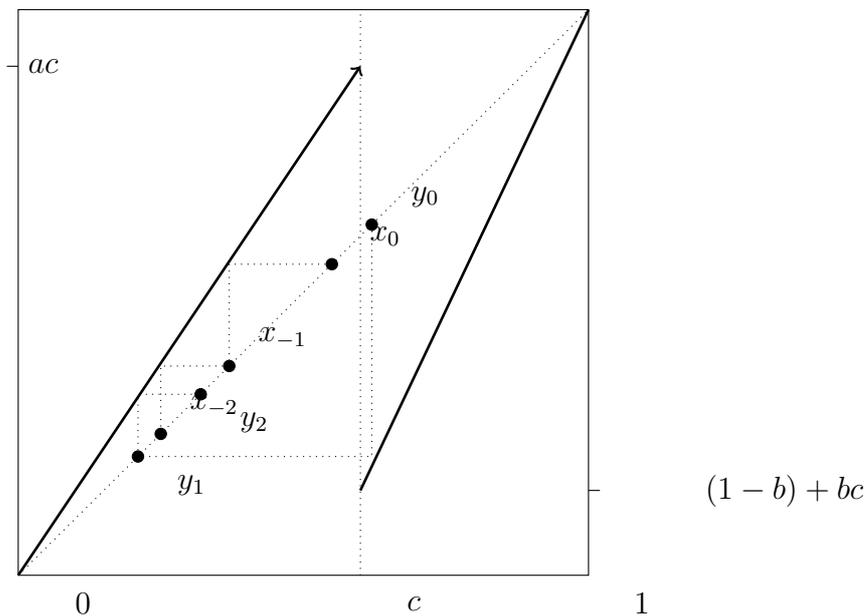

\begin{proposition} The map $Tx:=\function{ax &\mbox{if~} 0\le x<c \\
                                    bx + 1-b  &\mbox{if~}c\le x\le1}$ 
with $0<c<1$ satisfies the $G(\phi)$ property with $\sum_k\phi(k)<\infty$ 
if and only if $\min(a,b)>1$ and $ac=b(1-c)=1$. 
\end{proposition}

%\bcr{This result demonstrates that it is possible that $T\in (A+A)$ while $T\not\in (U+A)$. [??]}

\proof Assume that $\min(a,b)>1$ and $ac=b(1-c)=1$. Let $\v{x}, \v{y}$ be two 
arbitrary trajectories of this map. We define the gluing trajectory $\v{z}$ as follows: 
$z_0:=y_0, z_k:=T^ky_0$ for $k>0$, and $z_{k-1}:=T^{-1}_{x_{k-1}}z_k$ for $k\le0$ 
(see Fig.~\ref{f:gluing-ex}). 
In fact, since the map $T$ is piecewise expanding, there are no other options for $z_0$, 
similarly to the Example~\ref{e:hyp}.

By this construction $z_k=y_k ~\forall k\ge0$, while for negative $k$ the distances 
between $z_k$ and $x_k$ decrease at exponential rate, since each time we are applying for their 
calculations the same inverse branch of the expanding map $T$. 

Nevertheless, observe that the distance between $Tx_0$ and $Ty_0$ (i.e. the gap between 
the backward trajectory $\v{x}$ and the forward trajectory $\v{y}$) might be arbitrary 
close to $1$ when the points $x_0$ and $y_0$ are close to the point $c$. 

To demonstrate that if $ac<1$ or $b(1-c)<1$, while $\min(a,b)>1$, the gluing property breaks down, 
consider a pair of trajectories $\v{x}:=\{0\},~ \v{y}:=\{1\}$. Here by $\{0\}$ and $\{1\}$ we 
mean trajectories staying at fixed points $0$ and $1$ correspondingly. Since the map $T$ 
is expanding, there are no trajectories converging simultaneously in positive or negative time to these ones. 
Therefore, there can be no gluing. The important point here is that for each $x\in X$ 
both the left and right preimages under the action of $T$ are well defined.

It remains to consider the case $\min(a,b)\le1$. Without losing generality, assume that 
$a\le b$. Then under the assumption $\min(a,b)=1$, the interval $[0,c)$ is forward 
invariant, which implies that again the pair of trajectories $\v{x}:=\{0\},~ \v{y}:=\{1\}$ 
cannot be glued together. Now, if $\min(a,b)<1$ the same property holds true but 
due to a different reason: there is no trajectory converging to 1 when time goes to $\infty$ 
and to $0$ when time goes to $-\infty$.

Therefore if $\min(a,b))>1$ we may set  
$\phi(k):=\function{(\min(a,b))^{k} &\mbox{if~} k\le0 \\
                                                  0 &\mbox{if~} k\ge0 }$.   
\qed

So far in all examples the rate function $\phi$ has exponential tails. In general this is 
absolutely not the case, which is demonstrated in the following example. 

\begin{example}\label{ex:neut} (Neutral fixed points)
$X:=[0,1]$, $Tx:=\function{x + (1-c)(\frac{x}{c})^{1+\alpha} &\mbox{if } x\le c \\ 
                                  1 - T(\frac{c(1-x)}{1-c})  &\mbox{if } x>c}$. %$0<c,~0\le\alpha\le1$. 
\end{example} 

\begin{figure}\begin{center}
\begin{tikzpicture}[scale=0.5]
     \draw [-](0,0) to (10,0) to (10,10) to (0,10) to (0,0);  
     \node at (0.2,-0.5){$0$}; \node at (10,-0.5){1}; \node at (-0.3,9.7){1}; 
     \draw [dotted] (0,0) to (10,10); \draw [dotted] (6,0) to (6,10); \node at (6,-0.5){$c$};
     \draw [thick]  (0,0) .. controls (3,3.25) and (5,6) .. (6,10);  %Bezier
     \draw [thick]  (10,10) .. controls (8,7.75) and (7,6) .. (6,0);  %Bezier
\end{tikzpicture}\end{center}
\caption{Neutral mapping.}\label{f:neut-m} \end{figure}
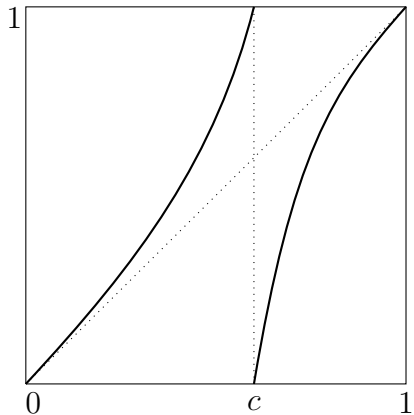

\begin{proposition}\label{p:neut} The map $T$ in the Example~\ref{ex:neut} 
(see Fig.~\ref{f:neut-m}) satisfies 
\begin{enumerate}
\item the (weak) gluing property (\ref{e:glu-w}) with the rate function 
        $\phi(k):=\function{C|k|^{-\gamma} &\mbox{if~} k\le0 \\
                                                    0 &\mbox{if~} k\ge0 }$, 
        where $C=C(\alpha,c)<\infty,~\gamma>1/\alpha$, and $\phi$ is summable if $0<\alpha<1$; 
\item if $\alpha>1$ there is no summable rate function $\phi$, for which $T\in G(\phi)$;
\item if $\alpha=0$ then $T\in G(\phi)$ with an exponentially decaying rate function.
\end{enumerate}
\end{proposition}

\begin{remark} The neutrally expanding map, considered in this example, for $0<\alpha<1$ 
satisfies only the weak gluing property (\ref{e:glu-w}), but the strong one (\ref{e:glu}) 
breaks down, which excludes the uniform (U+U) shadowing property.
\end{remark}

\proof The map $T$ is uniformly piecewise expanding everywhere except the neighborhoods 
of the neutral fixed points $0$ and $1$. Therefore the analysis of the forward part of the gluing 
property does not differ much from the situation in the example~\ref{e:plin},  
namely, we choose $z_0:=y_0$. Still we need to show that proper chosen pre-images 
of the point $z_0$ will approximate the backward semi-trajectory well enough. 

To this end we need some estimates of the rate of convergence in backward time of a map 
with a neutral fixed point.

\begin{lemma}\label{l:neut} Let $\tau(v):=v+Rv^{1+\alpha},~R>0, \alpha>0$. Then \\
(a) $\tau^{-n}(v) \le K n^{-\gamma}~~\forall v\in[0,1],~n\in\IZ_+$ 
     and some $K<\infty,~\gamma>1/\alpha$.\\
(b) $\tau^{-n}(v) \ge Kv n^{-\gamma}~~\forall v\in[0,1],~n\in\IZ_+$ 
     and some $K<\infty,~\gamma<1/\alpha$.\\
(c)  if $\alpha=0$ then $\tau^{-n}(v) \le (1+R)^{-n}v~~\forall v\in[0,1],~n\in\IZ_+$.\\
\end{lemma}
\proof First, we estimate $\tau^{-1}v$ using convexity of the function $\tau$. 

\begin{figure} \begin{center}
\begin{tikzpicture}[scale=0.75]
     \draw [-](0,0) to (10,0);  \draw [-](0,0) to (0,10);
     \draw [dotted] (0,0) to (10,10); \draw [dotted] (7,0) to (7,10); \draw [dotted] (7,7) to (0,7);
     \draw [-] (0,0) to (7,10); \draw [-] (4.5,4) to (7,10); \node at (7,-0.5){$v$}; \node at (7.7,9.8){$\tau(v)$}; 
     \draw [line width=1.5pt]  (0,0) .. controls (3,3.25) and (5,6) .. (7,10);  %Bezier
     \draw [dotted] (4.90,0) to (4.90,7); \node at (4.8,-0.5){$u$}; 
     \draw [dotted] (5.35,0) to (5.35,7); \node at (5.33,-0.5){$\t{v}$}; 
     \draw [dotted] (5.75,0) to (5.75,7); \node at (5.85,-0.5){$w$}; 
%     \draw [line width=1pt](0,0) to (6,9); 
\end{tikzpicture}\end{center}
\caption{Estimates for the inverse neutral mapping: 
             $u\le \t{v}:=\tau^{-1}v \le w \le v$.}\label{f:neut}\end{figure}
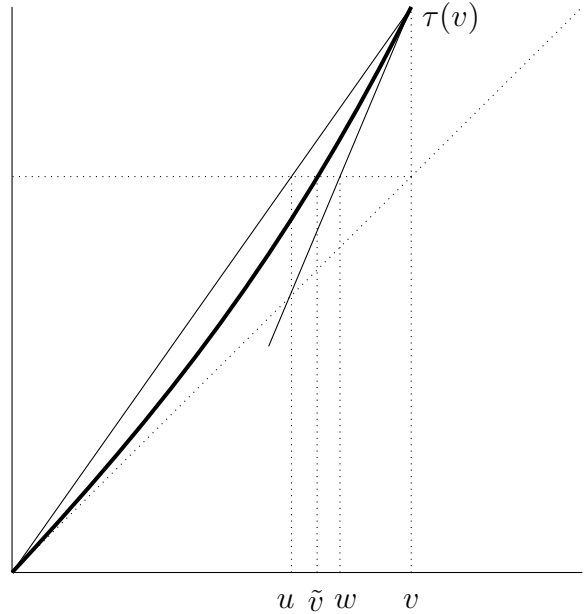
In Fig.~\ref{f:neut} the upper straight line connects the origin with the point $(v,\tau(v))$, 
while the lower straight line is tangent to the graph of the function $\tau$ at point $v$. 
A simple computation gives: %
\beq{e:est}{ u=\frac{v}{1+rv^\alpha} \le\tau^{-1}v 
       \le w=v\left( 1 - \frac{Rv^\alpha}{1 + (1+\alpha)Rv^{1+\alpha}} \right) .}%

Using (\ref{e:est}) we prove assertion (a) by induction.
$$ \tau^{-1}v \le v\left( 1 - \frac{Rv^\alpha}{1 + (1+\alpha)Rv^{1+\alpha}} \right) \le 1 $$ 
for $v\le 1$. Therefore, it is enough to set $K<1$. Assuming that 
$\tau^{-n}(v) \le K n^{-\gamma}$, we prove this for $n+1$: 
$$ \tau^{-n-1}(v) = \tau^{-1}(\tau^{-n}v) \le \tau^{-1}(K n^{-\gamma}) $$
$$ \le Kn^{-\gamma}\left(1 - \frac{R(Kn^{-\gamma})^\alpha}
              {1 + (1+\alpha)R(Kn^{-\gamma})^{1+\alpha}} \right) $$
$$ = K(n+1)^{-\gamma}(1+\frac1n)^{-\gamma}\left(1 - \frac{R(Kn^{-\gamma})^\alpha}
              {1 + (1+\alpha)R(Kn^{-\gamma})^{1+\alpha}} \right) .$$
The last expression is less or equal to $K(n+1)^{-\gamma}$, provided $\gamma\alpha>1$.

\bigskip

To prove assertion (b), we again use the inductive argument. 
$$ \tau^{-n-1}(v) = \tau^{-1}(\tau^{-n}v) \ge \tau^{-1}(K v n^{-\gamma}) 
     \ge \frac{Kvn^{-\gamma}}{1 + R(Kvn^{-\gamma})^\alpha} 
       = \frac{Kv(n+1)^{-\gamma} (1+\frac1n)^{-\gamma}}{1 + R(Kvn^{-\gamma})^\alpha} ,$$
which exceeds $Kv(n+1)^{-\gamma}$ if $\alpha\gamma>1$.  

\bigskip

Assertion (c) is a simple consequence the fact that $\tau'v=(1+R)>1$ in this case.
\qed

\n{\bf Continuation of the Proof of Proposition~\ref{p:neut}.} 
Applying the assertion (a) of Lemma~\ref{l:neut} to the inverse branches of the map $T$, 
and using that $0\le u\le1$ we get %
\beq{e:neut}{ \rho(T^{-n}u, \{0,1\}) < C n^{-\gamma}~~\forall u\in X,~n\in\IZ_+ ,} %
where $\rho(u,A):=\inf_{a\in A}\rho(u,a)$, $\rho(u,v):=|u-v|$, 
and $\gamma>1/\alpha,~ C=C(\alpha,c)<\infty$. 

Now we are ready to estimate $\rho(x_{-n},z_{-n})$ for $n\in \IZ_+$. 
By the triangle inequality, using (\ref{e:neut}), we get 
$$ \rho(x_{-n},z_{-n}) \le \rho(x_{-n},0) + \rho(0,z_{-n}) 
     \le 2C n^{-\gamma} .$$
Therefore, since $z_n\equiv y_n~\forall n\in\IZ_+$, 
$\phi(k):=\function{2C|k|^{-\gamma} &\mbox{if~} k\le0 \\
                                              0 &\mbox{if~} k\ge0 }$ 
defines the rate function for the (weak) gluing property (\ref{e:glu-w}).  
Moreover, if $\alpha<1$ then $\gamma>1$ and hence $\phi$ is summable. 

If $\alpha=0$, then it follows from Lemma~\ref{l:neut}(c) that $T\in G(\phi)$ 
with an exponentially decaying rate function.

Finally, the absence of the summable rate function $\phi$, for which $T\in G(\phi)$, 
follows from Lemma~\ref{l:neut}(b). \qed

%%%%%%%%%%%%%%%%%%%%%%%%%%%

\end{document}